\documentclass[12pt]{article}
\newcommand{\qed}{\hfill\rule{4pt}{8pt}\par\vspace{\baselineskip}}
\newtheorem{de}{Definition}[section]

\newtheorem{pr}[de]{Proposition}

\def\bea{\begin{eqnarray*}}\def\eea{\end{eqnarray*}}

\begin{document}
\title{Hopf bimodules are modules over a \\
diagonal crossed product algebra}

\author {Florin Panaite\\Institute of Mathematics of the Romanian Academy\\
P. O. Box 1-764, RO-70700 Bucharest, Romania\\e-mail: fpanaite@stoilow.imar.ro}

\date{}
\maketitle
\begin{abstract}If $H$ is a finite dimensional Hopf algebra, C. Cibils and 
M. Rosso found an algebra $X$ having the property that Hopf bimodules 
over $H^*$ coincide with left $X$-modules. We find two other algebras, 
$Y$ and $Z$, having the same property; namely, $Y$ is the $two-sided\;
crossed\;product$ $H^*\# (H\otimes H^{op})\# H^{* op}$ and $Z$ is the 
$diagonal\;crossed\;product$ $(H^*\otimes H^{* op})\bowtie (H\otimes H^{op})$ 
(both concepts are due to F. Hausser and F. Nill). We also find explicit 
isomorphisms between the algebras $X, Y, Z$.
\end{abstract} 
\section{Introduction}
${\;\;\;}$Let $H$ be a finite dimensional Hopf algebra. In \cite{cr} 
C. Cibils and M. Rosso introduced an algebra $X=(H^{op}\otimes H)
\underline{\otimes }(H^*\otimes H^{* op})$ having the property that  
Hopf bimodules over $H^*$ coincide with left $X$-modules. This algebra $X$ 
was further used in \cite{t} by R. Taillefer who proved that if $M$ and $N$ 
are (finite dimensional) Hopf bimodules over $H^*$ then the 
Gerstenhaber-Schack cohomology groups $H^*_{GS}(M, N)$ are isomorphic to 
$Ext^*_X(M, N)$. The multiplication of $X$ is a ``twist'' of the one 
of $(H^{op}\otimes H)\otimes (H^*\otimes H^{* op})$, 
and it was also proved in \cite{cr} that $X$   
is isomorphic to the direct tensor product between a Heisenberg double and 
(the opposite of) a Drinfel'd double. \\
${\;\;\;}$In this paper we introduce two algebras, $Y$ and $Z$, with the  
same property as $X$ (Hopf bimodules over $H^*$ coincide with left  
modules over $Y$ or $Z$) but which are more ``structured'' than $X$. Namely, 
$Y$ is the $two-sided\;crossed\;product$ (in the sense of F. Hausser and 
F. Nill \cite{hn}) $H^*\#(H\otimes H^{op})\#H^{* op}$, and $Z$ is the  
$diagonal$ $crossed\;product$ (also in the sense of Hausser and Nill)  
$(H^*\otimes H^{* op})\bowtie (H\otimes H^{op})$.  
We also write down explicit   
isomorphisms between the algebras $X, Y, Z$ having the property that if  
$M$ is a Hopf bimodule over $H^*$ then the actions of $X, Y, Z$ on $M$  
correspond via these isomorphisms. \\
${\;\;\;}$Let us mention that, among these three algebras, our favourite 
is $Z$, because the formulae for its multiplication and action on a Hopf 
bimodule look more elegant than in the other two cases. However, our 
approach relies on the structure of the algebra $Y$ and the use of 
``three corners'' Hopf modules. \\
${\;\;\;}$Also, let us note that, as a consequence of the Maschke-type 
theorem for diagonal crossed products (see \cite{HN}, Th. 8.2) it 
follows that if $H$ is semisimple and cosemisimple then $Z$ is semisimple, 
so we obtain a proof for the semisimplicity of $X$ independent on the 
isomorphism $X\simeq $$\cal H$$(H^*)\otimes D(H^*)^{op}$.     
\section{The algebras $X, Y, Z$}
${\;\;\;}$Throughout, $k$ will 
be a fixed field and all algebras, linear spaces etc. 
will be over $k$; unadorned $\otimes $ means 
$\otimes _k$. For coalgebras and Hopf algebras we shall use the framework of 
\cite{sw}. In particular for coalgebras we shall use $\Sigma $-notation: 
$\Delta (h)=\sum h_1\otimes h_2$, $(id \otimes \Delta )(\Delta (h))=
(\Delta \otimes id )(\Delta (h))=\sum h_1\otimes h_2\otimes h_3$ etc. \\
${\;\;\;}$In what follows, $H$ will be a finite dimensional Hopf algebra 
with antipode $S$. We start by stating some results which are either 
well-known (see for instance \cite{mon}) or easy to prove. \\
${\;\;\;}$If $A$ is a left $H$-module algebra, that is $A$ is an algebra and 
also a left $H$-module with action denoted by $h\otimes a\mapsto h\cdot a$ 
and such that $h\cdot (ab)=\sum (h_1\cdot a)(h_2\cdot b)$ and $h\cdot 1_A=
\varepsilon (h)1_A$ for all $h\in H$ and $a, b\in A$, the (left) smash 
product $A\# H$ is the algebra structure on $A\otimes H$ given by 
$$(a\# h)(b\# g)=\sum a(h_1\cdot b)\#h_2g $$
${\;\;\;}$A left $A\# H$-module may be identified with a vector space $M$ 
which is a left $H$-module and a left $A$-module (actions denoted by 
$h\otimes m\mapsto h\cdot m$ and $a\otimes m\mapsto a\cdot m$) related by the 
following compatibility condition:
$$h\cdot (a\cdot m)=\sum (h_1\cdot a)\cdot (h_2\cdot m)$$
for all $h\in H$, $a\in A$, $m\in M$.\\
${\;\;\;}$If $B$ is a right $H$-module algebra, that is $B$ is an algebra and 
also a right $H$-module with action denoted by $b\otimes h\mapsto b\cdot h$ 
and such that $(ab)\cdot h=\sum (a\cdot h_1)(b\cdot h_2)$ and $1_B\cdot h=
\varepsilon (h)1_B$ for all $h\in H$ and $a, b\in B$, the (right) smash 
product $H\# B$ is the algebra structure on $H\otimes B$ given by 
$$(h\# a)(g\# b)=\sum hg_1\#(a\cdot g_2)b$$
${\;\;\;}$We are not interested in right $H\# B$-modules, but also in left 
$H\# B$-modules. One can check that a left $H\# B$-module may be identified 
with a vector space $M$ which is a left $H$-module and a left $B$-module 
(actions denoted by $h\otimes m\mapsto h\cdot m$ and $b\otimes m\mapsto 
b\cdot m$) related by the following compatibility condition: 
$$b\cdot (h\cdot m)=\sum h_1\cdot ((b\cdot h_2)\cdot m)$$
for all $h\in H$, $b\in B$, $m\in M$.\\
${\;\;\;}$If $A$ is a right $H$-comodule algebra, that is $A$ is an algebra 
and a right $H$-comodule with structure map $\rho :A\rightarrow A\otimes H$ 
which is an algebra map, we can consider the categories $_A$$\cal M$$^H$ and 
$\cal M$$_A^H$ of relative Hopf modules (see \cite{mon}, p. 144). It is 
well-known (see \cite{doi}) that the category $_A$$\cal M$$^H$ may be 
identified with the category of left $A\# H^*$-modules ($A$ becomes a left 
$H^*$-module algebra as usual). Also (see \cite{mon}) $\cal M$$_A^H$ may be 
identified with the category of right $A\#H^*$-modules, but we need to 
identify it with a category of left modules, and this may be done as follows: 
since $A$ is a left $H^*$-module algebra, $A^{op}$ becomes a right 
$H^*$-module algebra with action given by $a\cdot p=S(p)\cdot a$ for all 
$a\in A$ and $p\in H^*$, where we denoted also by $S$ the antipode of $H^*$. 
Then one can prove that $\cal M$$_A^H$ may be identified with the category 
of left modules over the right smash product $H^*\#A^{op}$. \\
${\;\;\;}$We are interested in the categories of Hopf modules 
$_{H^*}^{H^*}$$\cal M$$^{H^*}$ and $^{H^*}$$\cal M$$_{H^*}^{H^*}$ (see for 
instance \cite{sch}) and mainly in the category of Hopf bimodules 
$_{H^*}^{H^*}$$\cal M$$_{H^*}^{H^*}$ (introduced for the first time in 
\cite{nic}). We have the obvious identifications:\\[2mm] 
$\;\;\;\;\;$$_{H^*}^{H^*}$$\cal M$
$^{H^*}\simeq _{H^*}$$\cal M$$^{H^*\otimes H^{* cop}}\simeq _{H^*}$
$\cal M$$^{(H\otimes H^{op})^*}$\\[2mm]
$\;\;\;\;\;$ $^{H^*}$$\cal M$$_{H^*}^{H^*}\simeq $
$\cal M$$_{H^*}^{H^*\otimes H^{* cop}}\simeq $$\cal M$$_{H^*}^{(H\otimes 
H^{op})^*}$\\[2mm]
where $H^*$ becomes a right $H^*\otimes H^{* cop}$-comodule   
algebra via the map $H^*\rightarrow H^*\otimes (H^*\otimes H^{* cop})$, 
$p\mapsto \sum p_2\otimes (p_3\otimes p_1)$. \\       
${\;\;\;}$Let $A$ be a left $H$-module algebra and $B$ a right $H$-module 
algebra (with actions denoted by $h\otimes a\mapsto h\cdot a$ and 
$b\otimes h\mapsto b\cdot h$). The $two-sided\;crossed\;product$ $A\# H\#B $,  
introduced by F. Hausser and F. Nill in \cite{hn}, is an algebra structure on 
$A\otimes H\otimes B$, given by 
$$(a\# h\# b)(a'\# h'\# b')=\sum a(h_1\cdot a')\# h_2h'_1\# (b\cdot h'_2)b'$$
(the unit is $1\# 1\# 1$). Obviously the natural maps from $H, A, B, A\# H, 
H\# B$ to $A\# H\# B$ are all algebra maps. \\
${\;\;\;}$Define the category $(A, H, B)-mod$ to be the category whose 
objects are vector spaces $M$ which are left $A$-modules, left $H$-modules 
and left $B$-modules, with actions denoted by $a\otimes m\mapsto a\cdot m$, 
$h\otimes m\mapsto h\cdot m$, $b\otimes m\mapsto b\cdot m$, related by the 
compatibility conditions: \\[2mm]
$(i)\;\;\;b\cdot (a\cdot m)=a\cdot (b\cdot m)$\\[2mm]
$(ii)\;\;\;b\cdot (h\cdot m)=\sum h_1\cdot ((b\cdot h_2)\cdot m)$\\[2mm]
$(iii)\;\;\;h\cdot (a\cdot m)=\sum (h_1\cdot a)\cdot (h_2\cdot m)$\\[2mm]
for all $a\in A, b\in B, h\in H, m\in M$. The morphisms are the maps which 
are $A$-linear, $H$-linear and $B$-linear. Let us also note that the 
conditions $(ii)$ and $(iii)$ above are respectively equivalent to 
$$h\cdot (b\cdot m)=\sum (b\cdot S^{-1}(h_2))\cdot (h_1\cdot m)$$
$$a\cdot (h\cdot m)=\sum h_2\cdot ((S^{-1}(h_1)\cdot a)\cdot m)$$
${\;\;\;}$Now we can describe the category of left $A\# H\# B$-modules. 
\begin{pr} There is a natural isomorphism of categories 
$$A\# H\# B-mod\simeq (A, H, B)-mod$$
\end{pr}
{\bf Proof:} The identifications are given as follows: 
$$a\cdot m=(a\# 1\# 1)\cdot m$$
$$h\cdot m=(1\# h\# 1)\cdot m$$
$$b\cdot m=(1\# 1\# b)\cdot m$$
and conversely 
$$(a\# h\# b)\cdot m=a\cdot (h\cdot (b\cdot m))$$
for all $a\in A, h\in H, b\in B, m\in M$. We shall only prove that the formula 
$(a\# h\# b)\cdot m=a\cdot (h\cdot (b\cdot m))$ gives indeed a left 
$A\# H\# B$-module structure on $M$ provided $(i)$, $(ii)$ and $(iii)$ are  
satisfied, and leave the rest to the reader. We calculate:\\[2mm]
$(a\# h\#b )\cdot ((a'\# h'\# b')\cdot m)=a\cdot (h\cdot (b\cdot (a'\cdot 
(h'\cdot (b'\cdot m)))))\\[2mm]
=a\cdot (h\cdot (a'\cdot (b\cdot (h'\cdot (b'\cdot m)))))$\\[2mm]
(using $(i)$)\\[2mm]
$=\sum a\cdot (h\cdot (a'\cdot (h'_1\cdot ((b\cdot h'_2)\cdot (b'\cdot m)))))$
\\[2mm]
(using $(ii)$)\\[2mm]
$=\sum a\cdot (h\cdot (a'\cdot (h'_1\cdot (((b\cdot h'_2)b')\cdot m))))\\[2mm]
=\sum a\cdot ((h_1\cdot a')\cdot (h_2h'_1\cdot (((b\cdot h'_2)b')\cdot m)))$
\\[2mm]
(using $(iii)$)\\[2mm]
$=\sum (a(h_1\cdot a'))\cdot (h_2h'_1\cdot (((b\cdot h'_2)b')\cdot m))\\[2mm]
=(\sum a(h_1\cdot a')\# h_2h'_1\# (b\cdot h'_2)b')\cdot m\\[2mm]
=((a\# h\# b)(a'\# h'\# b'))\cdot m$, q.e.d.\qed  

${\;\;\;}$From the above discussion it follows that a left $A\# H\# B$-module 
is a left $H$-module $M$ which is also a left $A$-module and a left  
$B$-module such that $a\cdot (b\cdot m)=b\cdot (a\cdot m)$ for all $a\in A, 
b\in B, m\in M$ and such that $M$ is also a left $A\# H$-module and a left 
$H\# B$-module. \\
${\;\;\;}$Define the algebra $Y=H^*\# (H\otimes H^{op})\# H^{* op}$, where 
$H^*$ is a left $H\otimes H^{op}$-module algebra with action
$$(h\otimes h')\cdot f=h \rightharpoonup f\leftharpoonup h'$$
for all $h, h'\in H$ and $f\in H^*$, where $\rightharpoonup $ and 
$\leftharpoonup $ are the regular actions of $H$ on $H^*$, given by 
$(h\rightharpoonup f)(h')=f(h'h)$, $(f\leftharpoonup h')(h)=f(h'h)$, and 
$H^{* op}$ is a right $H\otimes H^{op}$-module algebra with action 
$$f\cdot (h\otimes h')=S(h\otimes h')\cdot f=(S(h)\otimes S^{-1}(h'))\cdot f=
S(h)\rightharpoonup f\leftharpoonup S^{-1}(h')$$
${\;\;\;}$So, the multiplication in $Y$ is given by:
$$(p\# (h\otimes g)\# q)(p'\# (h'\otimes g')\# q')=$$
$$=\sum p(h_1\rightharpoonup 
p'\leftharpoonup g_1)\# (h_2h'_1\otimes g'_1g_2)\# q'(S(h'_2)\rightharpoonup 
q\leftharpoonup S^{-1}(g'_2))$$
where the multiplications on the last two positions are made in $H$ and 
$H^*$ ($not$ in $H^{op}$ and $H^{* op}$).\\
${\;\;\;}$Now we come to Hopf bimodules. It is clear that a $H^*$-Hopf 
bimodule is a left $H\otimes H^{op}$-module $M$ (i.e. a $H^*$-bicomodule)  
which is also an $H^*$-bimodule and such that $M$ is an object in the 
categories $_{H^*}^{H^*}$$\cal M$$^{H^*}$ and $^{H^*}$$\cal M$$_{H^*}^{H^*}$. 
Regarding the right $H^*$-module structure of $M$ as a left $H^{* op}$-module 
structure, it is clear that $a\cdot (b\cdot m)=b\cdot (a\cdot m)$ for all 
$a\in A=H^*$ and $b\in B=H^{* op}$. \\
${\;\;\;}$Now, from the identifications\\[2mm] 
$_{H^*}^{H^*}$$\cal M$$^{H^*}\simeq 
_{H^*}$$\cal M$$^{(H\otimes H^{op})^*}\simeq _{H^*\# (H\otimes H^{op})}$
$\cal M$\\[2mm]
$^{H^*}$$\cal M$$_{H^*}^{H^*}\simeq $$\cal M$$_{H^*}^{(H\otimes H^{op})^*}
\simeq _{(H\otimes H^{op})\#H^{* op}}$$\cal M$ \\[2mm]
and from all the above, we obtain finally:
\begin{pr} There is a natural isomorphism of categories between 
$_{H^*}^{H^*}$$\cal M$$_{H^*}^{H^*}$ and the category of left 
$H^* \# (H\otimes H^{op})\# H^{* op}$-modules. 
\end{pr}
${\;\;\;}$We write down explicitly the $Y$-module structure of a Hopf  
bimodule $M\in _{H^*}^{H^*}$$\cal M$$_{H^*}^{H^*}$. Denote by $p\otimes 
m\mapsto p\cdot m$ and $m\otimes q\mapsto m\cdot q$ the $H^*$-bimodule 
structure of $M$ and by $M\rightarrow H^*\otimes M\otimes H^*$, $m\mapsto 
\sum m_{(-1)}\otimes m_{(0)}\otimes m_{(1)}$ the $H^*$-bicomodule structure 
of $M$ (the left $H\otimes H^{op}$-module structure of $M$ is then given by 
$(h\otimes g)\cdot m=\sum m_{(-1)}(g)m_{(1)}(h)m_{(0)}$). Then the action 
of $Y$ on $M$ is given by:\\[2mm]
$(p\# (h\otimes g)\# q)\cdot m=p\cdot ((h\otimes g)\cdot (m\cdot q))\\[2mm]
=\sum p\cdot ((m\cdot q)_{(-1)}(g)(m\cdot q)_{(1)}(h)(m\cdot q)_{(0)})\\[2mm]
=\sum p\cdot ((m_{(-1)}q_1)(g)(m_{(1)}q_3)(h)m_{(0)}\cdot q_2)\\[2mm]
=\sum m_{(-1)}(g_1)q_1(g_2)m_{(1)}(h_1)q_3(h_2)p\cdot m_{(0)}\cdot q_2\\[2mm]
=\sum m_{(-1)}(g_1)m_{(1)}(h_1)p\cdot m_{(0)}\cdot (h_2\rightharpoonup q
\leftharpoonup g_2)$\\[2mm]
for all $p, q\in H^*$, $h, g\in H$ and $m\in M$.\\
${\;\;\;}$Recall now from \cite{cr}, \cite{t} the structure of the algebra 
$X$ of Cibils and Rosso (we write it for $H^*$, that is in the formulae 
in \cite{cr} and \cite{t} one has to take $A=H^*$). The algebra structure 
of $X$ is $X=(H^{op}\otimes H)\underline{\otimes }(H^*\otimes H^{* op})$, 
where the multiplication is defined such that the first two and last two 
tensorands keep natural multiplication, 
$(g\otimes h)\underline{\otimes }(p\otimes q)=((g\otimes h)
\underline{\otimes }(1\otimes 1))((1\otimes 1)\underline{\otimes }(p\otimes 
q))$ and 
$$((1\otimes 1)\underline{\otimes }(p\otimes q))((g\otimes h)
\underline{\otimes }(1\otimes 1))=$$
$$=\sum p_1(S(g_1))p_3(S^{-1}(h_1))q_1(S^{-1}(g_3))q_3(S(h_3))
((g_2\otimes h_2)\underline{\otimes }(p_2\otimes q_2))=$$
$$=\sum (g_2\otimes h_2)\underline{\otimes }(S^{-1}(h_1)\rightharpoonup p
\leftharpoonup S(g_1)\otimes S(h_3)\rightharpoonup q\leftharpoonup 
S^{-1}(g_3))$$
${\;\;\;}$If $M$ is a Hopf bimodule over $H^*$, with notation as above, 
$M$ becomes a left $X$-module with action given by:
$$((g\otimes h)\underline{\otimes }(p\otimes q))\cdot m=
\sum m_{(-1)}(g_2)m_{(1)}(h_2)(h_1\rightharpoonup p\leftharpoonup g_1)\cdot 
m_{(0)}\cdot (h_3\rightharpoonup q\leftharpoonup g_3)$$
${\;\;\;}$Now, if we look at this formula and the one of the action of $Y$ 
on $M$, it is quite clear how to define an algebra isomorphism between $X$ 
and $Y$, such that the actions on $M$ correspond. 
\begin{pr} The map $\varphi :X\rightarrow Y$, given by
$$\varphi ((g\otimes h)\underline{\otimes }(p\otimes q))=\sum h_1
\rightharpoonup p\leftharpoonup g_1\# (h_2\otimes g_2)\# q$$
is an algebra isomorphism, having the property that 
$$((g\otimes h)\underline{\otimes }(p\otimes q))\cdot m=\varphi ((g\otimes h) 
\underline{\otimes }(p\otimes q))\cdot m$$
for all $g, h\in H$, $p, q\in H^*$ and $m$ in a $H^*$-Hopf bimodule $M$. 
The inverse of $\varphi $ is given by $\varphi ^{-1}:Y\rightarrow X$, 
$$\varphi ^{-1}(p\# (h\otimes g)\# q)=\sum (g_2\otimes h_2)
\underline{\otimes }(S^{-1}(h_1)\rightharpoonup p\leftharpoonup S(g_1)
\otimes q)$$
\end{pr}
{\bf Proof:} We shall only prove that 
$$\varphi (((1\otimes 1)\underline{\otimes }(p\otimes q))((g\otimes h)
\underline{\otimes }(1\otimes 1)))=$$
$$=\varphi ((1\otimes 1)\underline{\otimes }(p\otimes q))\varphi  
((g\otimes h)\underline{\otimes }(1\otimes 1))$$
and leave the rest of the computations to the reader. We calculate:\\[2mm]
$\varphi (((1\otimes 1)\underline{\otimes }(p\otimes q))((g\otimes h)
\underline{\otimes }(1\otimes 1)))=\\[2mm]
=\varphi (\sum (g_2\otimes h_2)\underline{\otimes }(S^{-1}(h_1)
\rightharpoonup p\leftharpoonup S(g_1)\otimes S(h_3)\rightharpoonup q
\leftharpoonup S^{-1}(g_3)))\\[2mm]
=\sum (h_2)_1S^{-1}(h_1)\rightharpoonup p\leftharpoonup S(g_1)(g_2)_1\# 
((h_2)_2\otimes (g_2)_2)\#S(h_3)\rightharpoonup q\leftharpoonup 
S^{-1}(g_3)\\[2mm]
=\sum p\# (h_1\otimes g_1)\#S(h_2)\rightharpoonup q\leftharpoonup 
S^{-1}(g_2)\\[2mm]
=(p\# (1\otimes 1)\#q)(1\# (h\otimes g)\# 1)\\[2mm]
=\varphi ((1\otimes 1)\underline{\otimes }(p\otimes q))\varphi  
((g\otimes h)\underline{\otimes }(1\otimes 1))$, q.e.d.\qed 
${\;\;\;}$Recall from \cite{hn} the definition of the $diagonal\; 
crossed\;product$ (which, in a slightly different form, appears also in 
\cite{wl} under the name $right\;twisted\;smash\;product$). If $C$ is an 
$H$-bimodule algebra with actions denoted by $h\otimes c\mapsto h\cdot c$ 
and $c\otimes h\mapsto c\cdot h$, the diagonal crossed product is the 
following associative algebra structure on $C\otimes H$:
$$(c\otimes h)(c'\otimes h')=\sum c(h_1\cdot c'\cdot S^{-1}(h_3))
\otimes h_2h'$$
${\;\;\;}$This structure is denoted by $C\bowtie H$; its unit is $1\bowtie 1$ 
and it contains $C\equiv C\bowtie 1$ and $H\equiv 1\bowtie H$ as 
subalgebras. \\
${\;\;\;}$As noted in \cite{wl}, a linear space $M$ is a left 
$C\bowtie H$-module if and only if it is a left 
$H$-module and a left $C$-module with  
actions $h\otimes m\mapsto h\cdot m$ and $c\otimes m\mapsto c\cdot m$ such 
that $h\cdot (c\cdot m)=\sum (h_1\cdot c\cdot S^{-1}(h_3))\cdot 
(h_2\cdot m)$ (the $C\bowtie H$-module structure of $M$ is then given 
by $(c\bowtie h)\cdot m=c\cdot (h\cdot m)$).\\
${\;\;\;}$If $A$ is a left $H$-module algebra and $B$ is a right $H$-module 
algebra, it was proved in \cite{hn} that $C=A\otimes B$ is an $H$-bimodule 
algebra with actions $h\cdot (a\otimes b)\cdot g=h\cdot a\otimes b\cdot g$ 
for all $a\in A, b\in B, h, g\in H$, and the map $f:A\# H\# B\rightarrow  
(A\otimes B)\bowtie H$ given by 
$$f(a\# h\# b)=((a\otimes 1)\bowtie h)((1\otimes b)\bowtie 1)=\sum 
(a\otimes b\cdot S^{-1}(h_2))\bowtie h_1$$
is an algebra isomorphism, with inverse $f^{-1}:(A\otimes B)\bowtie H
\rightarrow A\# H\# B$, given by 
$$f^{-1}((a\otimes b)\bowtie h)=(1\# 1\# b)(a\# h\# 1)=\sum a\# h_1\# 
b\cdot h_2$$                
${\;\;\;}$From the description of modules over two-sided crossed products 
and over diagonal crossed products one can see that left modules over   
$A\# H\# B$ coincide with the ones over $(A\otimes B)\bowtie H$ and that 
$f(a\# h\# b)\cdot m=(a\# h\# b)\cdot m$ for all $a\in A, b\in B, h\in H$ and 
$m\in M$, where $M$ is such a module. \\
${\;\;\;}$We have seen that $H^*$ is a left $H\otimes H^{op}$-module algebra 
and $H^{* op}$ is a right $H\otimes H^{op}$-module algebra, so we can 
consider the diagonal crossed product $Z=(H^*\otimes H^{* op})\bowtie 
(H\otimes H^{op})$, whose multiplication may be written as:
$$((p\otimes q)\bowtie (h\otimes g))((p'\otimes q')\bowtie (h'\otimes g'))=$$
$$=\sum (p(h_1\rightharpoonup p'\leftharpoonup g_1)\otimes 
(h_3\rightharpoonup q'\leftharpoonup g_3)q)\bowtie (h_2h'\otimes g'g_2)$$
(where, as above, the products in the second and fourth positions are 
in $H^*$ and $H$, $not$ in $H^{* op}$ and $H^{op}$).\\
${\;\;\;}$From the above discussion, the modules over $Z$ are also the same 
as Hopf bi\-mo\-dules over $H^*$, and the algebras $Y$ and $Z$ are isomorphic  
via the maps $\alpha :Y\rightarrow Z$, $\alpha ^{-1}:Z\rightarrow Y$, 
$$\alpha (p\# (h\otimes g)\# q)=\sum (p\otimes (h_2\rightharpoonup q
\leftharpoonup g_2))\bowtie (h_1\otimes g_1)$$
$$\alpha ^{-1}((p\otimes q)\bowtie (h\otimes g))=\sum p\# (h_1\otimes g_1)\# 
S(h_2)\rightharpoonup q\leftharpoonup S^{-1}(g_2)$$
${\;\;\;}$Hence, the algebras $X$ and $Z$ are also isomorphic, via the maps 
$\beta :X\rightarrow Z$, $\beta ^{-1}:Z\rightarrow X$, 
$\beta =\alpha \circ \varphi $, 
$\beta ^{-1}=\varphi ^{-1}\circ \alpha ^{-1}$, that is 
$$\beta ((g\otimes h)\underline{\otimes }(p\otimes q))=
\sum (h_1\rightharpoonup p\leftharpoonup g_1\otimes h_3\rightharpoonup q
\leftharpoonup g_3)\bowtie (h_2\otimes g_2)$$
$$\beta ^{-1}((p\otimes q)\bowtie (h\otimes g))=\sum (g_2\otimes h_2)
\underline{\otimes }(S^{-1}(h_1)\rightharpoonup p\leftharpoonup S(g_1)
\otimes S(h_3)\rightharpoonup q\leftharpoonup S^{-1}(g_3))$$
and, if $M$ is an $H^*$-Hopf bimodule, the actions of $X$ and $Z$ on $M$ 
correspond via these isomorphisms. The action of $Z$ on $M$ is given by 
$$((p\otimes q)\bowtie (h\otimes g))\cdot m=\sum m_{(-1)}(g)m_{(1)}(h)
p\cdot m_{(0)}\cdot q$$        

\end{document}